\documentclass[twoside,12pt]{amsart}
\date{\today}
\usepackage{enumerate}
\usepackage{graphicx}
\usepackage{amsmath}
\usepackage{verbatim}  % I added this (bw)
\usepackage{amssymb}
\usepackage{epsfig}
\usepackage{amsthm}%entornos theorem
\usepackage{amscd}
\usepackage[usenames]{color}

\definecolor{r}{rgb}{.6,0,.3}

\newcounter{intro}

\usepackage[latin1]{inputenc}

\newcommand{\D}{\mathbb{D}}

\newcommand{\C}{\mathbb{C}}

\newcommand{\R}{\mathbb{R}}
\newcommand{\RR}{\mathbb{R}}

\newcommand{\h}{\mathbb{H}}
\newcommand{\hs}{\mathfrak{H}}
\newcommand{\fD}{\mathfrak{D}}

\newcommand{\ka}{\kappa}
\newcommand{\esf}{\mathbb{S}}

\newcommand{\calC}{\mathcal{C}}
\newcommand{\calP}{\mathcal{P}}

\renewcommand{\Im}{\operatorname{Im}}

\newcommand{\whC}{\widehat{C}}

\newcommand{\del}{\partial}

\newtheorem{theorem}{Theorem}[section]

\newtheorem{proposition}[theorem]{Proposition}
\newtheorem*{maintheorem*}{Main Theorem}
\newtheorem*{theorem*}{Theorem}
\newtheorem*{mainlemma*}{Main Lemma}
\newtheorem*{keylemma*}{Key Lemma}

\theoremstyle{definition}
\newtheorem{remark}[theorem]{Remark}
\newtheorem*{remark*}{Remark}

\title[Families of minimal surfaces in $\h^2 \times \R$ foliated by arcs]
{Families of minimal surfaces in $\h^2 \times \R$ foliated by
arcs and their Jacobi fields}%
\author[L. Ferrer]{Leonor Ferrer}
\address[Ferrer]{
  Departamento de Geometr\'\i{}a y Topolog\'\i{}a,
  Universidad de Granada,
  18071 Granada, Spain.
}
\email{lferrer@ugr.es}
\author[F. Martín]{Francisco Mart\'\i{}n}
\address[Martín]{
  Departamento de Geometr\'\i{}a y Topolog\'\i{}a,
  Universidad de Granada,
  18071 Granada, Spain.
}
\email{fmartin@ugr.es}
\author[R. Mazzeo]{Rafe Mazzeo}
\address[Mazzeo]{%  
Department of Mathematics,
Stanford University,
Stanford, California 94305, U.S.A.
}
\email{mazzeo@math.stanford.edu}
\author[M. Rodríguez]{Magdalena Rodríguez}
\address[Rodríguez]{
  Departamento de Geometr\'\i{}a y Topolog\'\i{}a,
  Universidad de Granada,
  18071 Granada, Spain.
}
\email{magdarp@ugr.es}
\thanks{
 L. Ferrer, F. Martín and M.M. Rodríguez are partially supported by the
  MINECO/FEDER grant MTM2014-52368-P and MTM2017-89677-P;  R. Mazzeo supported by
the NSF grant DMS-1608223; F. Martin is also partially supported by the
Leverhulme Trust grant IN-2016-019.
    }
\begin{document}

\maketitle
\begin{abstract}
  This note provides some new perspectives and calculations regarding an interesting known family of minimal  surfaces
  in $\h^2 \times \R$.  The surfaces in this family are the catenoids, parabolic catenoids and tall rectangles.
  Each is foliated by either circles, horocycles or circular arcs in horizontal copies of $\h^2$.  All of these surfaces
  are well-known, but the emphasis here is on their unifying features and the fact that they lie in a single continuous
  family. We also initiate a study of the Jacobi operator on the parabolic catenoid, and compute the Jacobi fields associated
  to deformations to either of the two other types of surfaces in this family.  
\end{abstract}

\section{Introduction}
In these notes we study properties of an interesting family of minimal surfaces in $\h^2 \times \RR$. The surfaces in this
family are foliated by circles or circular arcs in parallel slices of $\h^2 \times \RR$; those surfaces foliated by entire circles are
catenoids, those foliated by horocycles are  parabolic catenoids, and
those foliated by arcs equidistant from a geodesic are the so-called tall rectangles.  Somewhat surprisingly, using the
Poincaré disk model of $\h^2$, these surfaces all appear as intersections of regular surfaces in $\R^3$ with the
unit cylinder $\{ (x,y,t) \in \RR^3 \; : \; x^2+y^2 < 1\}$. 

\section{Preliminaries and notation}
To set notation, we shall use both the Poincar\'e disk and upper half-space models of $\h^2$; these
are denoted by $\D=\{ z \in \C \ :\ |z|<1\}$ and $\hs=\{z=x+{\rm i}\, y \in\C  \,: \; y>0\}$. We also denote by
$o=(0,0) \in \D$.
The isometry between these two models is the M\"obius transformation $g:\D \rightarrow \hs$,
$$
g(z)={\rm i } \; \frac{1-z}{1+z}.
$$
These have metrics 
\[
d\rho^2 = 4 \, \frac{|dz|^2}{(1-|z|^2)^2}, \quad \mbox{and} \qquad \frac{1}{y^2} \; |dz|^2,
\]
respectively; the corresponding metric for $\h^2 \times \R$ is $d\sigma^2 = d\rho^2 + dt^2$.
%\[
%d\sigma^2=4 \, \frac{|dz|^2}{(1-|z|^2)^2}+d\, t^2= d \, \rho^2+d\, t^2\, .
%\]
We denote by $\del \h^2$ the boundary at infinity of $\h^2$, usually identified
with the boundary $S^1=\{ z \in \C \ :\ |z|=1\}$  of $\D$, but sometimes also with the boundary $\R$ of $\hs$ together with the extra point
at infinity. Recall from \cite{KM} that there are several interesting compactifications of the boundary at
infinity of $\h^2 \times \R$, but for the purposes of this paper we consider only the portion $\del \h^2 \times \R$, which
lies in the boundary of any of these useful compactifications.

\section{Catenoids of revolution in $\h^2 \times \R$} \label{cat}
We begin by studying the minimal surfaces of revolution in $\h^2 \times \RR$. These analogues of
catenoids in $\R^3$ were originally described in the seminal article \cite{n-r} of Nelli and Rosenberg.
We take a slightly different approach to their construction here. Consider a conformal harmonic
parametrization of a minimal annulus $\mathbb A$:  
$$
X = (F, h): \Delta_R  =\{R< |z|<1\} \longrightarrow \mathbb A \subset \h^2 \times \R;
$$
thus
$$
F:\Delta_R \to \h^2, \qquad h : \Delta_R \to \R
$$
are both harmonic maps. We impose the condition that $h$ is locally constant on $\del \Delta_R$, i.e., takes two
different constant values on the two boundary components of this annulus; by translation we assume that
these two values are $0$ and $h_0 > 0$.

We now pass to the induced mappings from the universal cover $M=\{ w=w_1+{\rm i} w_2 \in \C \; : \; 0 < w_2  < 1 \}$ of
$\Delta_R$. The (holomorphic) covering map is
$$
\varphi: M \longrightarrow \Delta_R, \qquad \varphi(w)=\exp \left(- {\rm i} \,(\log R) \,w\right).
$$
Now $\hat h:= h \circ \varphi: M \to \R$ is harmonic, and we normalize by assuming that
$\hat h|_{w_2=0} \equiv 0$, $\hat h|_{w_2=1} \equiv h_0$. Since it is bounded in the lateral directions, $\hat{h}(w) = h_0 \Im w$;
this is the imaginary part of the holomorphic function $z = h_0 w$ which, following notation in \cite{n-r}, we write as $z= -\theta + it$.
This is defined on the strip $M_{h_0}=\{ w=w_1+{\rm i} w_2 \in \C \; : \; 0 < w_2  < h_0 \}$. 
% and $\Omega=z(M)=\{ \zeta \in \C \, : \; 0<{\rm Im}(\zeta)<h_0 \}$.

%Assuming that $\mathbb A$ is transverse to any horizontal slice $\h^2 \times \{t_0\}$, $z$ is a local conformal
%coordinate.
We now bring in the fact that $\mathbb A$ is a surface of revolution, so we can write
\begin{equation}
\label{eq:F}
\hat F(\theta, t):=(F \circ \varphi \circ z^{-1})(\theta, t) =r(t) \, e^{i s(\theta)}
\end{equation}
for some smooth functions $s$ and $r$. The Hopf differential equals $\hat \Phi= \frac 14 dz^2$. Since
the pair $(\hat F, \hat{h}) = (\hat{F}, t)$ is conformal on $M_{h_0}$, we also have that
$$
\| \hat F_t \|^2_{\h^2}+1= \| \hat F_\theta \|^2_{\h^2},
$$
which becomes
$$
\frac{4}{(1-r(t)^2)^2} \left( s'(\theta)^2r(t)^2-r'(t)^2\right)=1.
$$
Rearranging this we obtain
$$
s'(\theta)^2= \frac{(1-r(t)^2)^2}{4 r(t)^2}+\left( \frac{r'(t)}{r(t)} \right)^2,
$$
hence each of the two sides must equal the same constant $\kappa^2$. Write $s'(\theta)=\kappa > 0$ and 
\begin{equation} \label{eq:C}
\kappa^2= \frac{(1-r(t)^2)^2}{4 r(t)^2}+\left( \frac{r'(t)}{r(t)} \right)^2.
\end{equation}
The harmonicity of $\hat F$ yields the equation
\begin{equation} \label{jod}
  \hat F_{z \, \overline{z}}+2  (\log \rho \circ \hat F)_u \hat F_z \, \hat F_{\overline{z}}=0,
\end{equation}
where $u$ is a (holomorphic) coordinate of $\h^2$ and the metric $d\rho^2$ on $\h^2$ equals $\rho(u)^2 |du|^2$.
We compute that 
$$
2 (\log \rho \circ \hat F)_u = 2 \overline{\hat F}\big/(1-| \hat F |^2),
$$
and hence the left hand side of \eqref{eq:C} becomes
\begin{multline}
%\hat F_{z \, \overline{z}}+2 (\log \rho \circ \hat F)_u \hat F_z \, \hat F_{\overline{z}}=  \\ 
\frac{1}{4} e^{{\rm i} s(\theta )} \left(r''(t)-r(t) \left(s'(\theta )^2-i
   s''(\theta )\right)\right)+ \\ \frac{r(t) e^{i s(\theta )} \left(r(t)^2
   s'(\theta )^2-r'(t)^2\right)}{2 \left(r(t)^2-1\right)}
\end{multline}
Substituting from \eqref{eq:C} and using that $s''(\theta)=0$, we arrive at the expression
$$
\hat F_{z \, \overline{z}}+2 (\log \rho \circ \hat F)_u \hat F_z \, \hat F_{\overline{z}}=  \frac{e^{{\rm i} s(\theta )}
  \left(4 r(t) r''(t)-4 r'(t)^2+r(t)^4-1\right)}{16    r(t)}.
$$
We conclude finally that \eqref{jod} is equivalent to
\begin{equation} \label{jod1}
4 r(t) r''(t)-4 r'(t)^2+r(t)^4-1=0.
\end{equation}
This is the equation obtained in a different way by Nelli and Rosenberg in \cite{n-r}.
\begin{figure}[htbp]
\begin{center}
\includegraphics[height=.66\textheight]{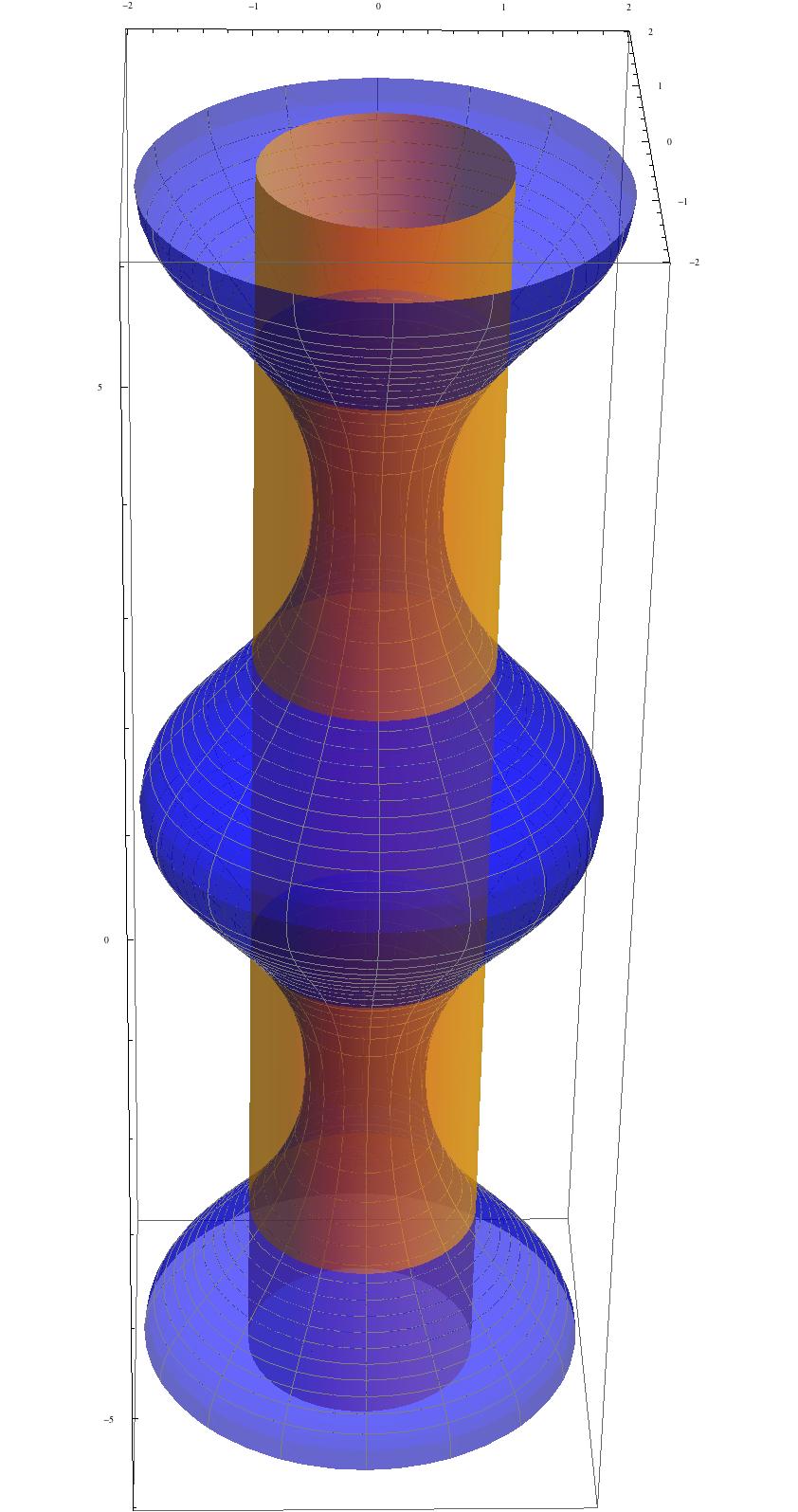}
\end{center}
\caption{One of the unduloid-type surfaces $U_\kappa$ whose intersection with the solid cylinder $\{ (x,y,t) \in \RR^3 \; : \; x^2+y^2 < 1\}$ gives the catenoids
in $\h^2 \times \R$.} \label{fig:one}
\end{figure}

From the fact that \eqref{eq:C} is a first integral of \eqref{jod1}, it is easy to read off that solutions $r_\kappa(t)$
are defined on the entire real line and are periodic in $t$, oscillating between two values:
$$
\sqrt{\kappa^2+1}-\kappa \leq r_\kappa(t) \leq \sqrt{\kappa^2+1}+\kappa.
$$
Noting that $\sqrt{\kappa^2+1}-\kappa< 1 <\sqrt{\kappa^2+1}+\kappa$, we see that $(\hat{F}, \hat{h})$
extends to a map from $\C$ to $\R^3$, with image a complete surface of revolution which we write as $U_\kappa$;
this has the conformal type of $\C^*$. This family of surfaces is similar in many ways to the classical family of
Delaunay unduloids. 

\begin{proposition}
The surfaces $U_\kappa$ converge to the cylinder $\{|z| = 1\} \times \R$, as $\kappa \to 0$, and to a foliation of $\R^3$ by
parallel planes as $\kappa \to \infty$. The connected components of $U_\kappa \cap \{ |z| <1\}$, which are all identified
with one another by appropriate vertical translations, are copies of the standard catenoid of revolution, which we
write as $\mathcal{C}_\kappa$, in $\h^2 \times \R$, see Figure \ref{fig:one}.  The surface $\mathcal{C}_\kappa$ is conformally
equivalent to a proper annulus $\Delta_R$ where $R = R_\kappa \in ({\rm e}^{-\pi}, 1)$. The height $h = h(\kappa)$
of $\mathcal C_\kappa$ decreases monotonically from $\pi$ to $0$ as $\kappa$ increases from $0$ to $\infty$. 
\end{proposition}

\section{Parabolic  Catenoids}
We next consider a family of surfaces obtained via a particular limit of horizontal dilations of the catenoids $\mathcal C_\kappa$. 

For any point $p \in \h^2$ and $h \in (0,\pi)$, let $\calC_{h,p}$ denote the catenoid in $\h^2 \times \RR$ which
is rotationally symmetric around the axis $\{p\} \times \RR$, symmetric with respect to reflections across $t=0$
and has height $h \in (0,\pi)$.  Observe that $\calC_{h(\kappa),p}$ is obtained applying an horizontal dilation to $\mathcal C_\kappa$.

\begin{figure}[htpb]
{\includegraphics[width=4cm]{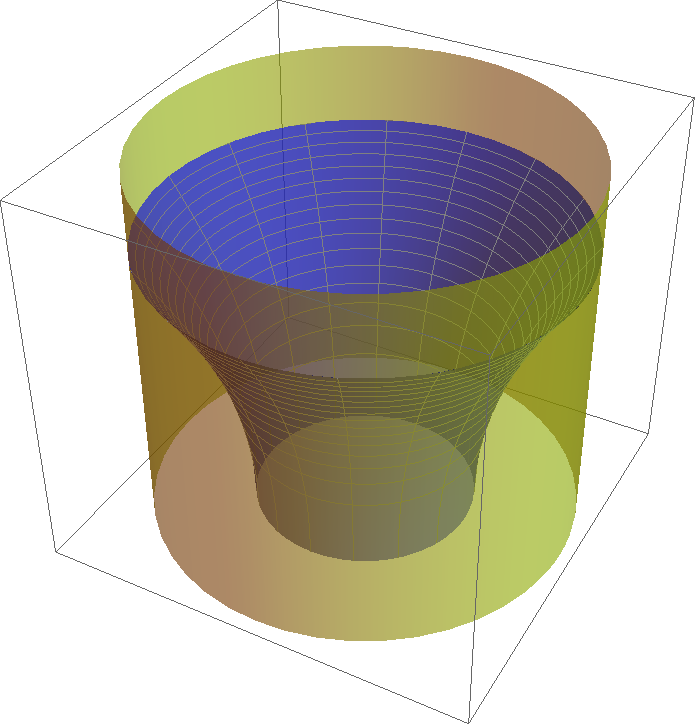}}
{\includegraphics[width=4cm]{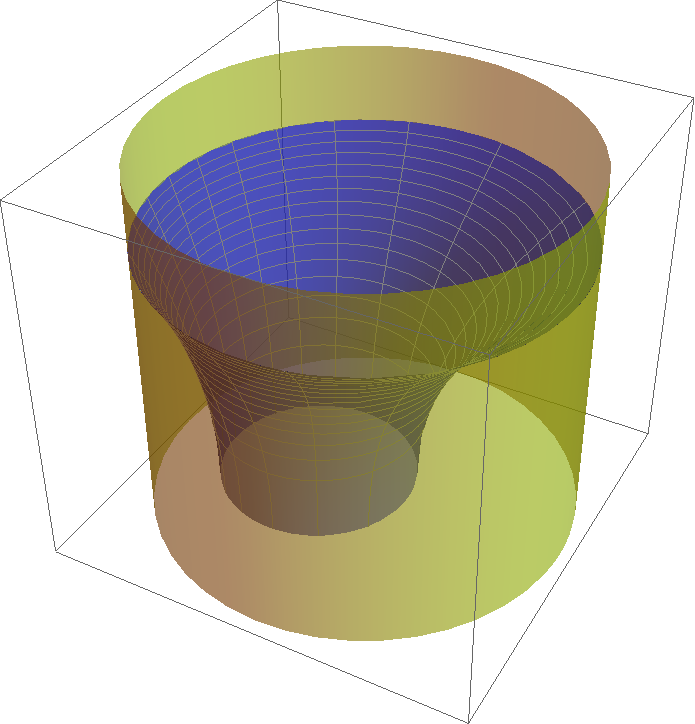}}
{\includegraphics[width=4cm]{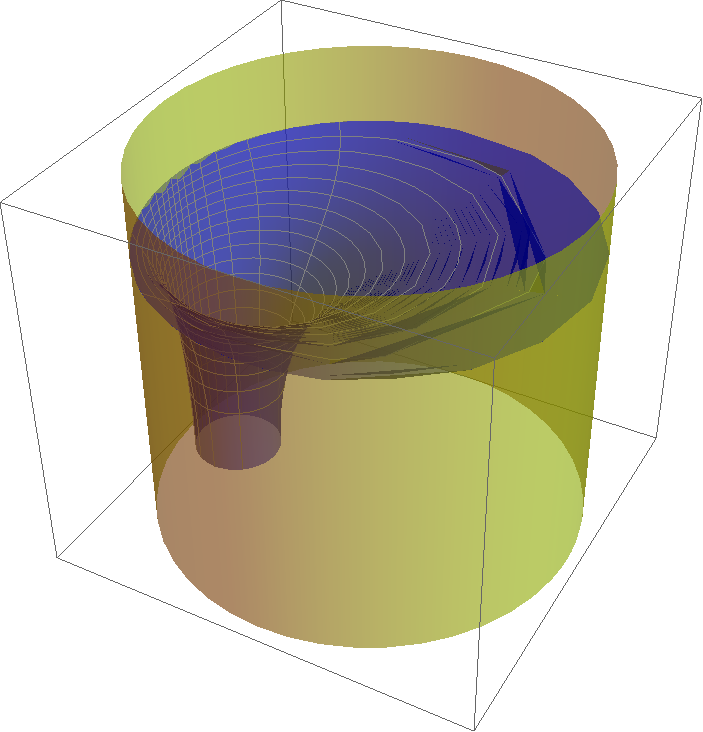}}
\caption{The upper half of three surfaces in the sequence $\mathcal{C}_j$. The limit is the union of two horizontal disks.}
\label{sequence}
\end{figure}

Now take a sequence of these
catenoids, $\calC_j := \calC_{h_j, p_j}$ such that $p_j \to q \in \del \h^2$ and $h_j$ remains bounded away from both $0$ and $\pi$.
Then $\calC_j$ converges locally in $\mathcal C^\infty$ on any compact set of $\h^2 \times \R$ to two horizontal disks
$\h^2 \times \{t_j\}$, $j = 1, 2$, where $|t_2 - t_1| = \lim h_j$, see Figure \ref{sequence}. (If $h_j \to 0$, then $\calC_j$ converges to one horizontal disk with multiplicity two.)  

Suppose however that we let $h_j \nearrow \pi$. Depending on the rates at which $p_j \to q$ and
$h_j \to \pi$, various possibilities can occur. We suppose that these sequences are balanced in such a way
that $\calC_j$ intersects a fixed compact set $K \subset \h^2 \times \R$ for every $j$; this is easy to arrange by an elementary
argument. In this case, standard results imply that some subsequence of the $\calC_j$ converges locally in $\calC^\infty$
to a complete properly embedded minimal surface which we write as $\fD$.  By definition, $\fD$ is a parabolic
catenoid. It has asymptotic boundary equal to the union of two horizontal circles separated by distance $\pi$, together with a
vertical segment joining these two circles. This class of surfaces was discovered independently by Hauswirth \cite{h} and Daniel \cite{d}.

We can finesse the construction of the sequence $\calC_j$: start with a sequence of catenoids $\calC_{\kappa_j}$,
rotationally symmetric around the axis $\{o\} \times \R$, with $\kappa_j \searrow 0$, and lying
in the slab $|t| < \frac12 h_j$.  Let $\gamma$ denote the hyperbolic geodesic through $o$ and converging to $q \in
\del \h^2$, and let $\sigma_j$ denote a hyperbolic dilation toward $q$ along this geodesic which has
the property that, extending $\sigma_j$ to an isometry of $\h^2 \times \R$ which acts only on the first factor,
$\sigma_j(\calC_{\kappa_j})$ is tangent to the axis $\{o\} \times \R$ at the point $(o, 0)$ 
(i.e., on the central circle of the catenoid), with the neck of the catenoid lying between this point and $q$.
\begin{figure}[h]
{\includegraphics[width=3.8cm]{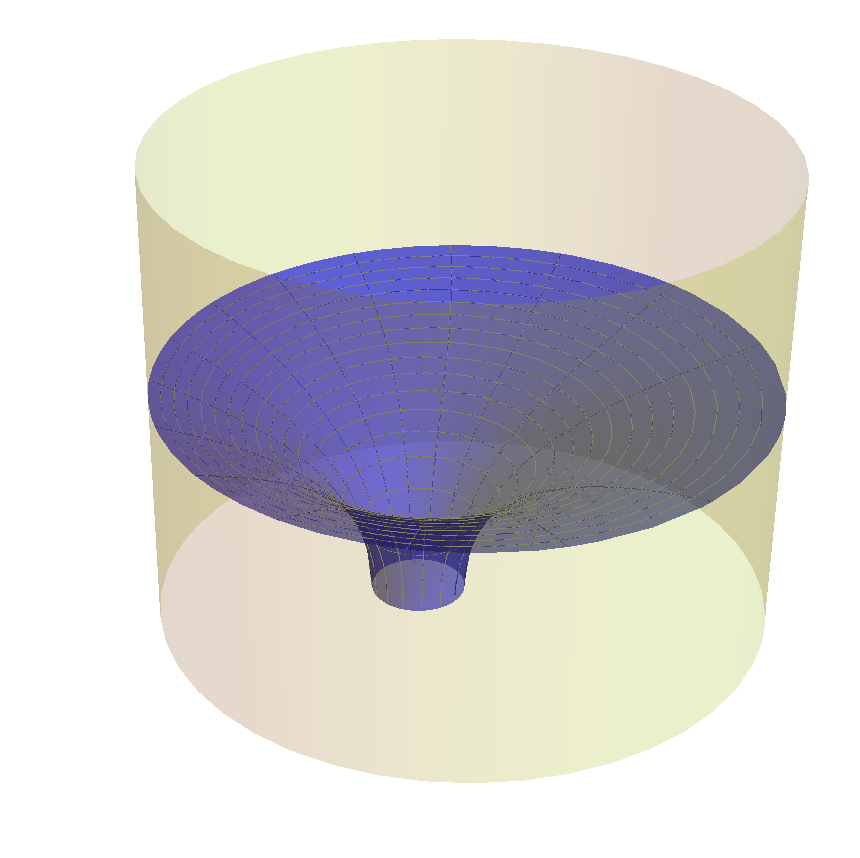}}
{\includegraphics[width=4cm]{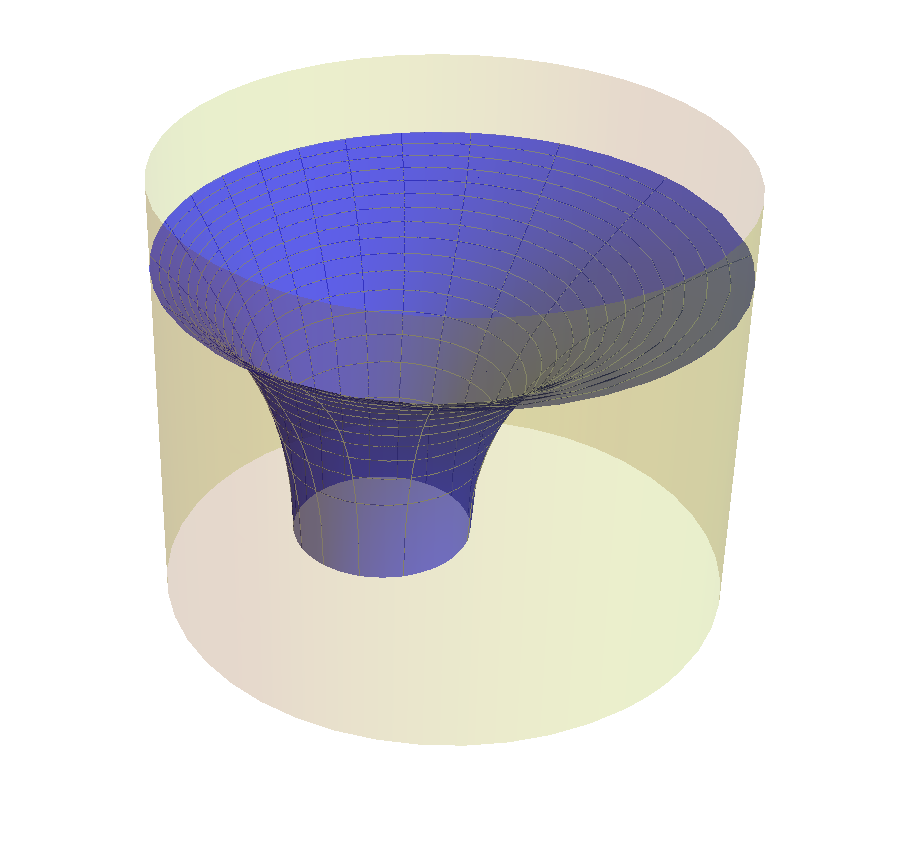}}
{\includegraphics[width=3.8cm]{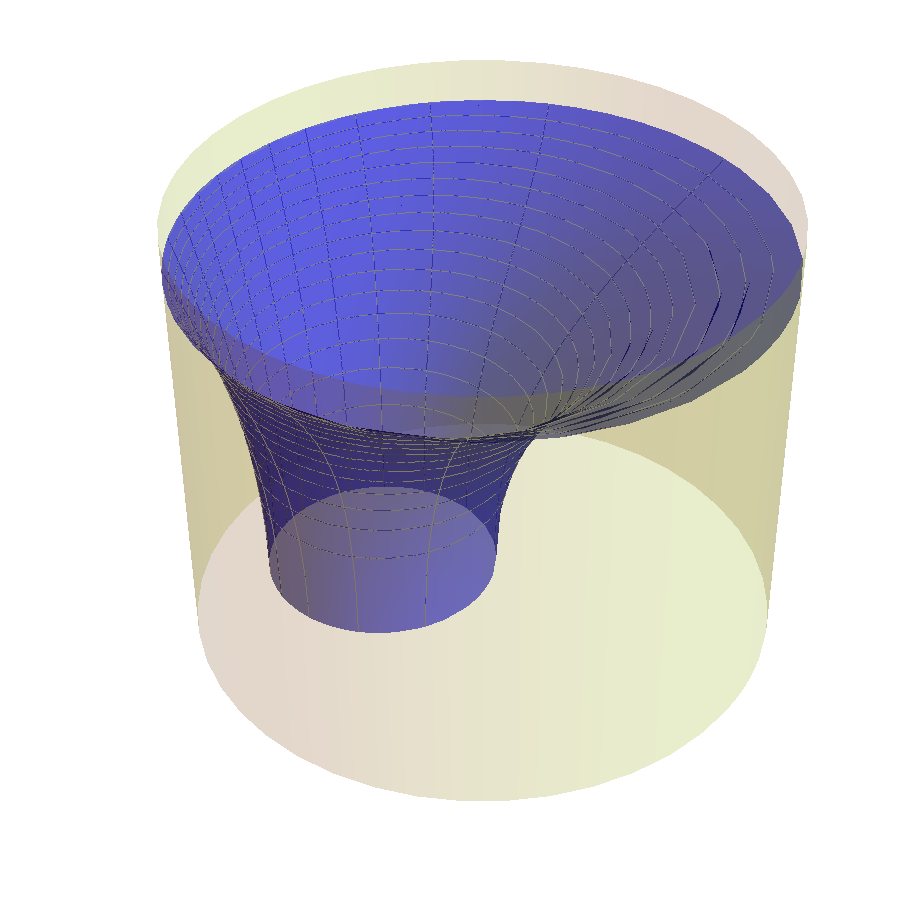}}
{\includegraphics[width=4cm]{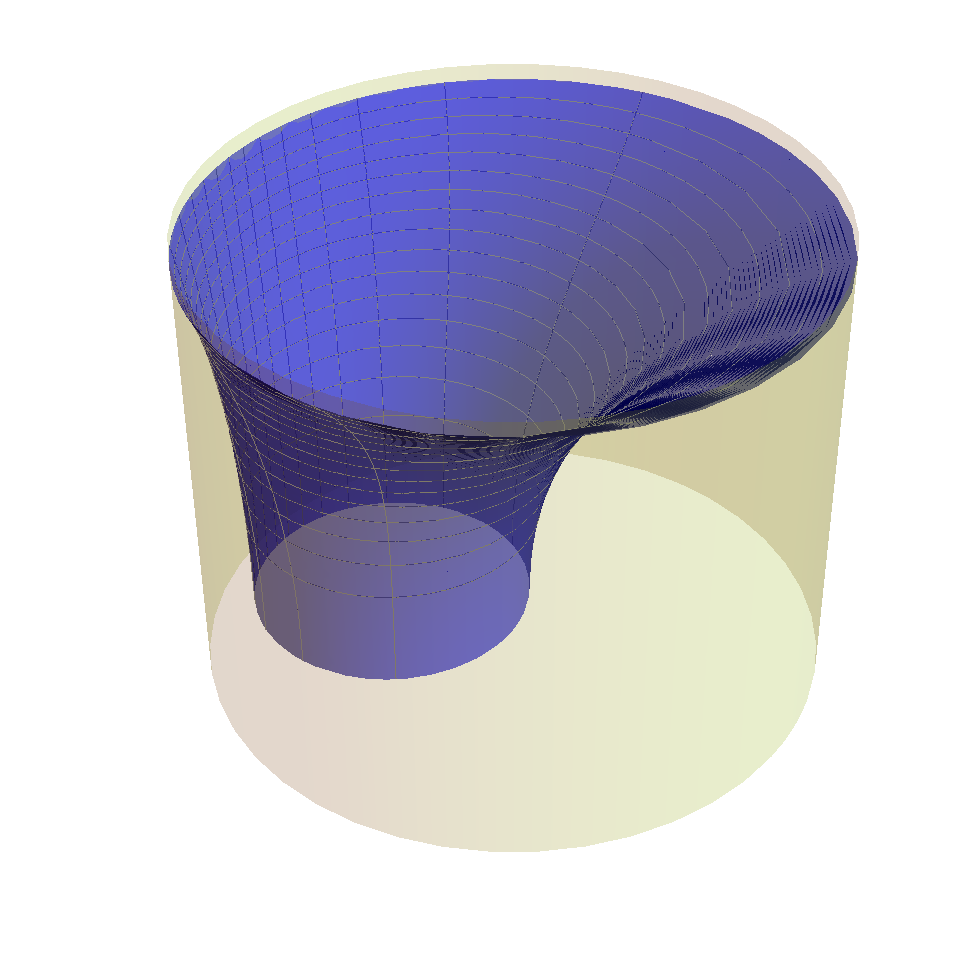}}
\caption{The upper half of four surfaces in the sequence $\sigma_j \left(\mathcal{C}_{\kappa_j}\right)$. The limit
  as $\kappa_j \searrow 0$ is the parabolic catenoid.}
\label{sequence2}
\end{figure}

This sequence certainly intersects a fixed compact set for all $j$, hence converges to a parabolic catenoid
which we denote by $\fD_q$,  see Figure \ref{sequence2}. Its height is precisely equal to $\pi$. The symmetries of the catenoid and the naturality
of the construction easily imply that these various Daniel surfaces are all related to one another by rotations of $\h^2
\times \R$ around the axis $\{o\} \times \R$.  We can also apply horizontal hyperbolic dilations to these surfaces,
which is the same as choosing a slightly different normalization of the sequence of surfaces $\calC_j$ above;
denoting the dilation parameter by $\lambda$, we obtain the family of surfaces $\fD_{q, \lambda}$; $\lambda = 1$
corresponds to the identity dilation. The surfaces in this entire family are all mutually isometric by rotations and hyperbolic dilations. 
We take as the standard model the surface in this family with $\lambda = 1$. It is an embedded disk with asymptotic
boundary the two horizontal circles $S^1 \times \{0\} \sqcup S^1 \times \{\pi\}$ and the vertical segment $\{q\} \times [0,\pi]$.

\begin{figure}[htpb]
{\includegraphics[width=6cm]{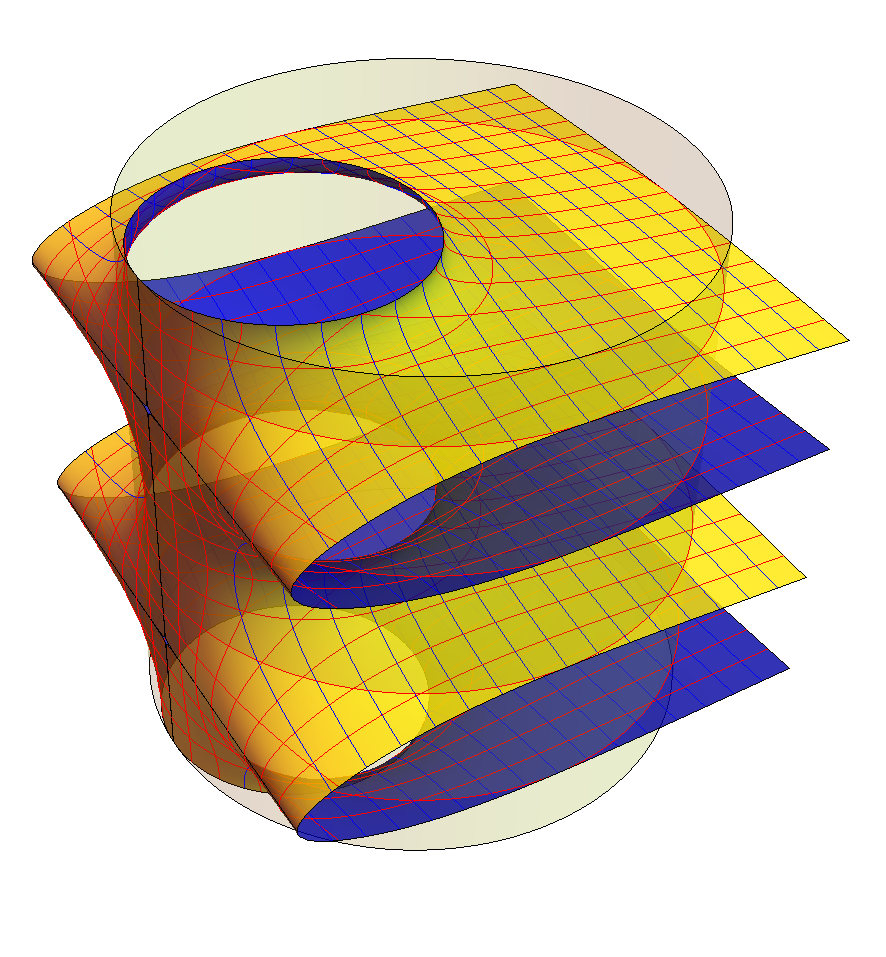}}
{\includegraphics[width=6cm]{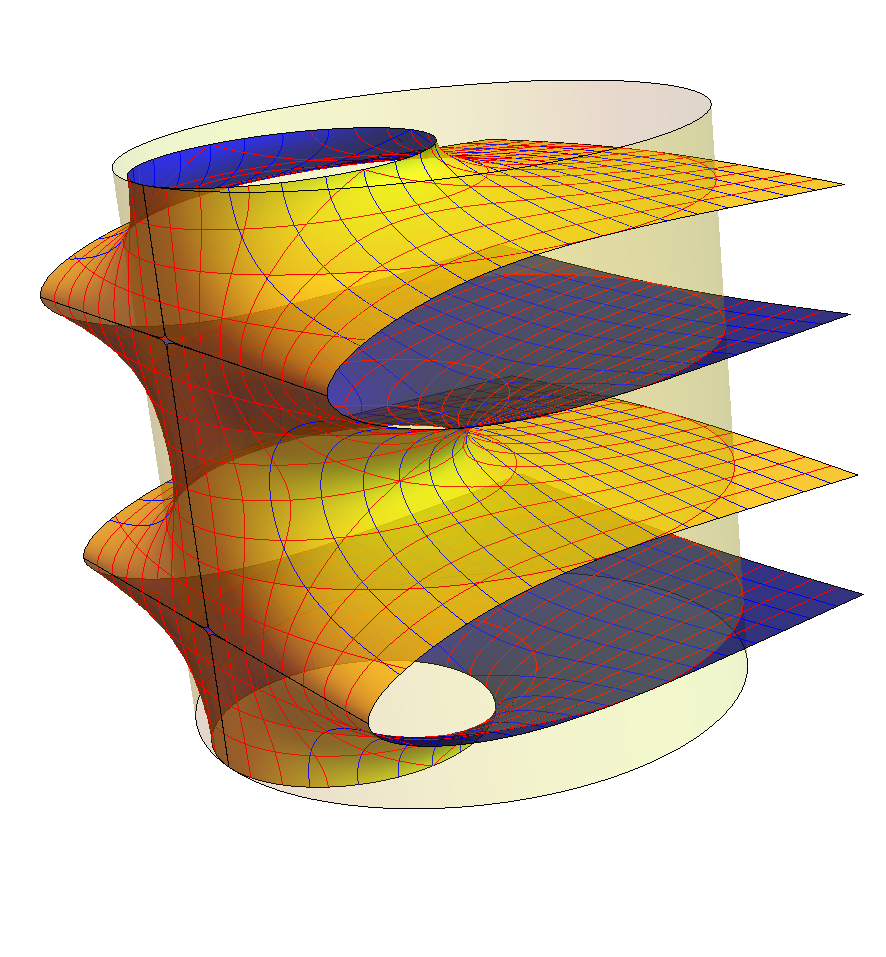}}
\caption{Two different views of the surface in $\R^3$ whose intersection with the solid cylinder $\{ (x,y,t) \in \RR^3 \; : \; x^2+y^2 < 1\}$ is the parabolic catenoid.} \label{fig:dani}
\end{figure}
\begin{remark}
  If one applies the limit process describe above to the entire unduloid-type surfaces $U_\kappa$, one obtains a complete
  limiting surface $\mathcal{Q}$ in $\R^3$. This limit (see Figure \ref{fig:dani}) is determined by
  $$
  \mathcal{Q}=\left\{(x,y,t) \in \R^3 \,: \, 1-x^2-y^2=\cos(t) ((1+x)^2+y^2)\right\}.
  $$
  The components of the intersection 
  $$\mathcal{Q} \cap \left\{(x,y,t) \in \R^3 \,: \, 1-x^2-y^2>0\right\}
  $$
  are an infinite union of (translated) copies of the parabolic catenoid.
\end{remark}

The simplest analytic representation of $\fD_{q,\lambda}$ uses the upper half-space model  for $\h^2$. Place $q$ at infinity in $\hs$
and use the standard coordinates $(x,y,t)$, $\hs \times \R$, where $y > 0$. The parabolic catenoid is invariant under parabolic translations,
which in this representation take the form $(x,y,t) \mapsto (x+b, y, t)$ for any $b \in \RR$, and it also intersects each horizontal slice $t = \mathrm{const.}$
transversely. Therefore this surface is a sweep-out of some curve $(0, f(t), t)$, $0 < t < \pi$ by these parabolic translations.  Searching for a minimal surface
with these properties leads in a straightforward way to the family of solutions $f(t) = \lambda \sin t$ for any $\lambda \in \RR^+$.
The surface is the image of corresponding family of embeddings of the strip $M_\pi = \R \times (0, \pi)$ given by
\begin{equation} \label{daniel}
\Psi_\lambda(x,t) = (\lambda \, x, \lambda \sin t, t). 
\end{equation}
For simplicity, we  write $\Psi = \Psi_1$.

\section{Tall Rectangles}
\begin{figure}[htbp]
\begin{center}
\includegraphics[height=.52\textheight]{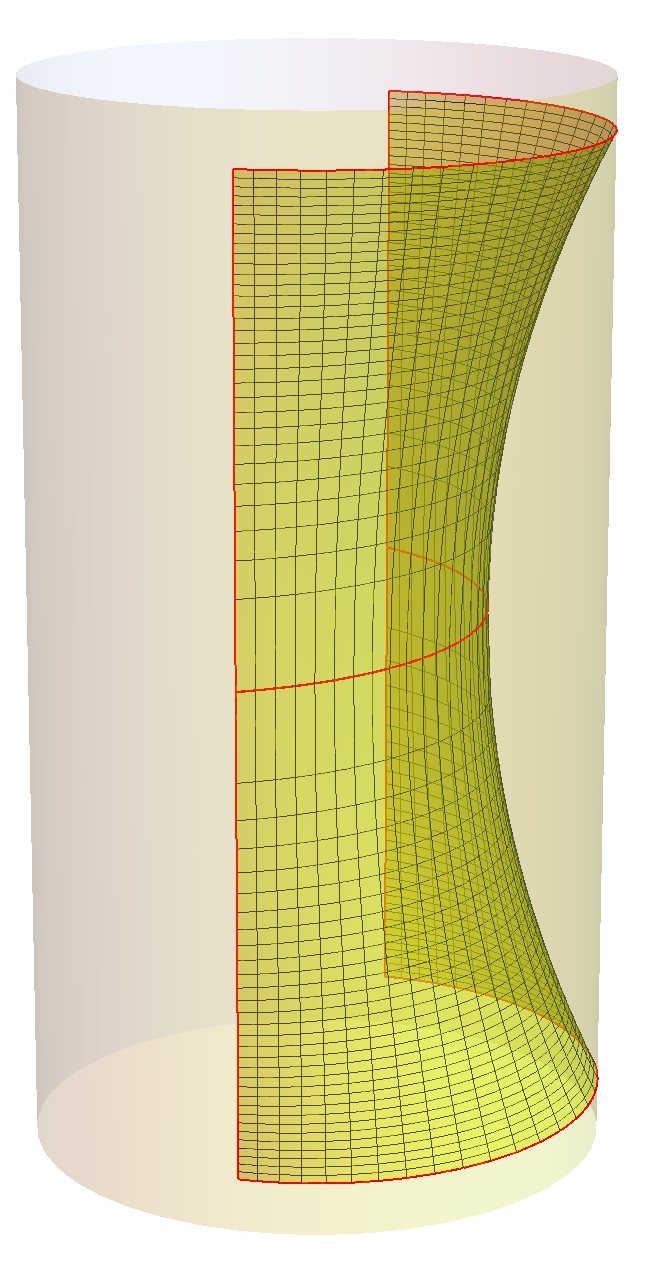}
\caption{The tall rectangle of height $h=\frac95 \pi$.}
\label{graph-11}
\end{center}
\end{figure}

The final family of surfaces in $\h^2 \times \R$ we consider here is the family of properly embedded mimimal disks, described
by Sa-Earp and Toubiana in \cite{st}. These have ideal boundary consisting of two parallel arcs $\sigma \times \{\pm h/2\}$,
where $h > \pi$ is arbitrary and $\sigma$ is an arc in $\partial \h^2$ with endpoints $q_1$ and $q_2$,
together with the vertical segments $\{ q_1\} \times [-h/2, h/2]$ and $\{q_2\} \times [-h/2, h/2]$, see Figure \ref{graph-11}.
Each of these surfaces, denoted $\Sigma_{\sigma, h}$, is area minimizing.  (A complete surface is area minimizing if any compact
piece is area-minimizing among all the surfaces with the same boundary.) The intersections $\Sigma_{\sigma, h} \cap( \h^2 \times \{t\})$
foliate $\Sigma_{\sigma, h}$ by a family of curves which, if all projected down to $\h^2$, all have the same endpoints $q_1$ and $q_2$ 
and are equidistant to the geodesic $[q_1,q_2]$. 

Following \cite{st}, we construct these surfaces as follows. Using the Poincaré disk model, consider the
vertical plane $\calP =\gamma \times \R$ where $\gamma$ is the geodesic $\{ \Im z = 0\} \subset\D$, parametrized
either by $x \in (-1,1)$ or by signed geodesic distance $\rho$ from $\{o\}$.
The relationship between the two parameters is $x=\tanh(\rho/2)$, $\rho=\log\left( \frac{1+x}{1-x}\right)$.
Next, fix $0 < d < 1$ and consider the curve $\sigma_d \subset \calP$ given as a bigraph of the two functions $t = \pm \lambda_d(\rho)$, i.e., 
$$
\sigma_d=\left\{(\rho, \pm \lambda_d(\rho)): \; \rho \geq \cosh^{-1}(1/d) \right\}.
$$
We then determine the conditions under which the surface swept out by this curve with respect to horizontal hyperbolic dilations
along the geodesic from $-i$ to $i$ is minimal.  A standard computations shows that minimality is equivalent to 
\begin{equation}
% \frac{d \lambda_e}{d\rho} (\rho)= \frac{e}{\sqrt{\cosh^2 \rho -e^2}}, \ \mbox{or}\ \
\frac{d\lambda_d}{d\rho}(\rho) = \frac{1}{\sqrt{d^2 \,\cosh^2 \rho -1}}
\end{equation}
(N.B.\ the treatment in \cite{st} uses $e = 1/d$ as a parameter instead.) By the chain rule, 
 \begin{equation}
\frac{d\lambda_d}{dx} = \frac{2}{\sqrt{d^2(1+x^2)^2-(1-x^2)^2}}.
 \end{equation}
The lower bound $\rho>\cosh^{-1}(1/d)$ transforms to $x>d_1 :=\frac{\sqrt{1-d}}{\sqrt{1+d}}$.
Note that the derivative is infinite at $x = d_1$, so this graph together with its reflection
cross the $x$-axis is at least $\calC^1$. A closer analysis shows that this bigraph is in fact $\calC^\infty$. 

In summary, the complete properly embedded minimal disk $\Sigma_d$ is the surface swept out by
the curve $\sigma_d$, where
\begin{equation}
  \lambda_d(x):= \int_{d_1}^{x} \frac{2 \; dv}{\sqrt{d^2(1+v^2)^2-(1-v^2)^2}}.
  \label{lambdad}
\end{equation}
It is straightforward to check that:
$$
\lambda_d(x)=-\frac{2}{1-d} \Im\left(F\left(\arcsin\left(d_1 x\right) \left| \frac{1}{d_1^2} \right)\right.\right),
$$
where %$F(\phi | z)$ is the classical elliptic integral of the first kind, given by:
$$
F(\phi | z):=\int_0^\phi \left(1- z \sin^2(\theta)\right)^{-1/2} \, d \theta, \quad -\frac \pi2 <\phi < \frac \pi2
$$
is the classical elliptic integral of the first kind. The surface itself is parametrized as a bigraph
$$
\Upsilon_d : (d_1,1) \times (-1,1) \rightarrow \h^2 \times \R,
$$
\begin{equation}
  \Upsilon_d(x, y)= \left(\frac{x-x y^2}{x^2 y^2+1},\frac{\left(x^2+1\right) y}{x^2 y^2+1}, \pm \lambda_d(x)\right).
  \label{Upsilon}
\end{equation}
%is a parametrization of $\Sigma_d^+=\Sigma_d \cap \{t>0\}$. Similarly, 
%$\Sigma_d^-=\Sigma_d \cap \{t<0\}$ can be parametrized by
%$$ \Upsilon_d^-(x, y)=
%\left(\frac{x-x y^2}{x^2 y^2+1},\frac{\left(x^2+1\right) y}{x^2 y^2+1}, -\lambda_d(x)\right) .$$
Note that the top and bottom halves $\Sigma_d^\pm = \Sigma_d \cap \{ \pm t \geq 0\}$
are each graphs over the ``lunette" region $L_d \subset \D$ lying between the circular
arc $\gamma_d$ passing through $\pm i$ and $d_1$, and the arc on the circumference joining
$-i$ to $i$ and passing through $1$. 

% $$
% \lambda_d(x):= \int_{d_1}^{\log( (x+1)/(x-1))} \frac{2 \; dv}{\sqrt{d^2(v^2+1)^2-(v^2-1)^2}}.
% $$
% (This is equivalent to the expression in \cite{st} after the change of variable $\rho = \log( (1+x)/(1-x))$. )
%Using classical properties of classical elliptic functions we have that
%\begin{multline}\lambda_d(x) =\\ -2 \Im\left(\frac{d \sqrt{1-\frac{(d-1) x^2}{d+1}} \sqrt{1-\frac{(d+1) x^2}{d-1}}
%   F\left(\sin ^{-1}\left(\sqrt{\frac{1-d}{1+d}}
%   x\right)|\frac{(d+1)^2}{(d-1)^2}\right)}{\sqrt{\frac{1-d}{1+d}} \sqrt{(d-1) x^2-d-1}
%   \sqrt{d \left(x^2-1\right)+x^2+1}}\right)\end{multline}
%   where $F(\phi | m)$ is the elliptic integral of the first kind given by
%   $$F(\phi | m)= \int_0^\phi (1-m \sin^2(\theta))^{-1/2} d \, \theta, \quad \mbox{ for $-\pi/2 <\phi <\pi/2.$} $$
% Next, joint this curve with the reflection of this graph across the $x$-axis. A short calculation shows that this
% bigraph is smooth, even at the intersection with the $x$-axis.   Note also that

The curve $\sigma_d$ and surface $\Sigma_d$ are contained in the region where
$$
|t| < \frac12 h_d = \int_{d_1}^1 \frac{2 \; dv}{\sqrt{d^2(v^2+1)^2-(v^2-1)^2}} < \infty.
$$
%The minimal surface $\Sigma_d^+$ is then obtained by applying all hyperbolic dilations along
%the geodesic connecting $-i$ to $i$ in $\D$.

% parametrized by translating the curve $\sigma_d$
% along the geodesic $[-{\rm i},{\rm i}]$ (which corresponds to the $y$-axis in the Poincaré
% disk.) If we reflect the above disk with respect to the plane $t=0$, then we obtain
% another disk  $\Sigma_d^-$ (which is also a vertical graph) and the Sa-Earp and
% Toubiana example is just the union of this two disks $\Sigma_d:=\Sigma_d^- \cup \Sigma_d^+$
% (see Figure \ref{graph-1}.)
It is not hard from this expression to check that the height $h_d$ of $\Sigma_d$ increases monotonically in $d$ and with
the following asymptotic behavior: 
\begin{itemize}
\item As $d \to 0$, $h_d \searrow \pi$ and $\Sigma_d$ diverges to infinity. Denoting by $\widetilde{T}_d$
the horizontal dilation along $\gamma$ which maps the point $d_1$ to $0$, then the family of disks
$Y_d:=(\widetilde{T}_{d} \times {\rm Id}_{\R})(\Sigma_d)$ converges to the parabolic catenoid passing through the origin.
\item As $d \to 1$, $h_d \to \infty$ and $\Sigma_d$ limits to the vertical plane $(-{\rm i}, {\rm i}) \times \R.$
\end{itemize}
\begin{figure}[h]
\begin{center}
\includegraphics[height=.52\textheight]{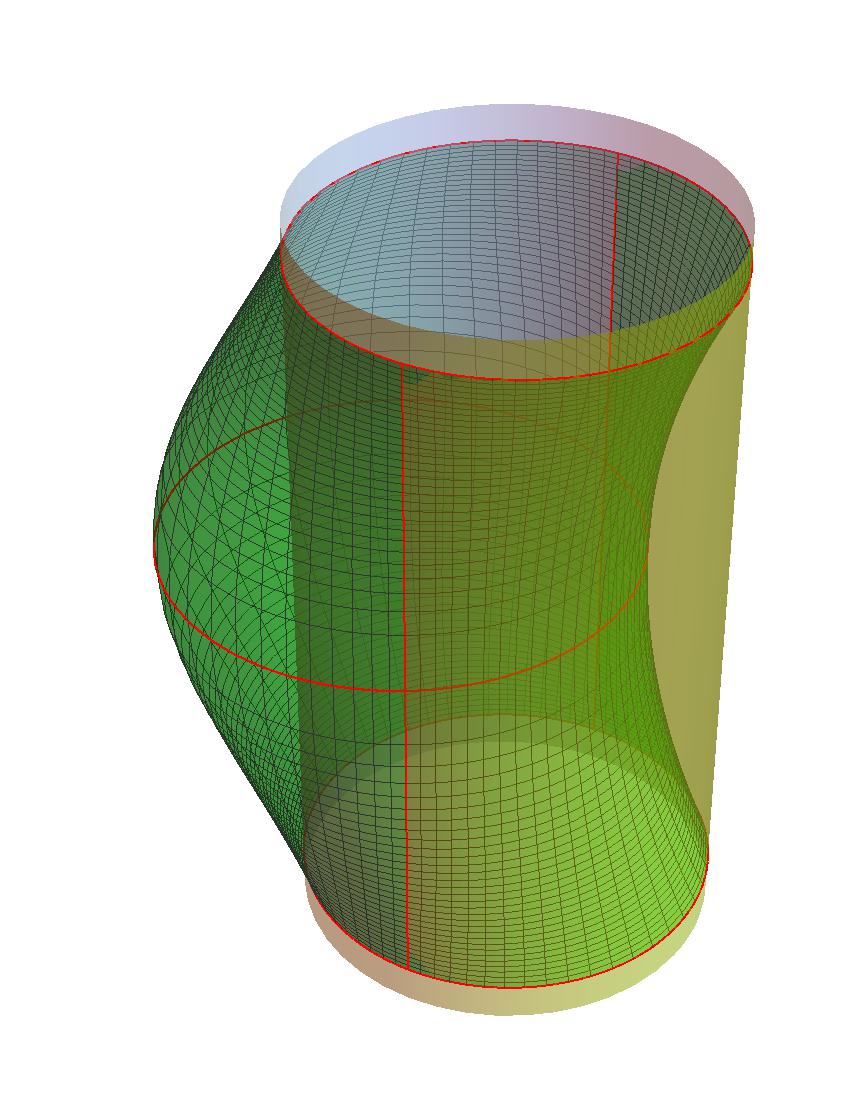}
\caption{The annulus  $\Upsilon_d((d_1,1/d_1) \times (-1,1))$.}
\label{annulus1}
\end{center}
\end{figure}

\begin{figure}[h]
\begin{center}
\includegraphics[height=.52\textheight]{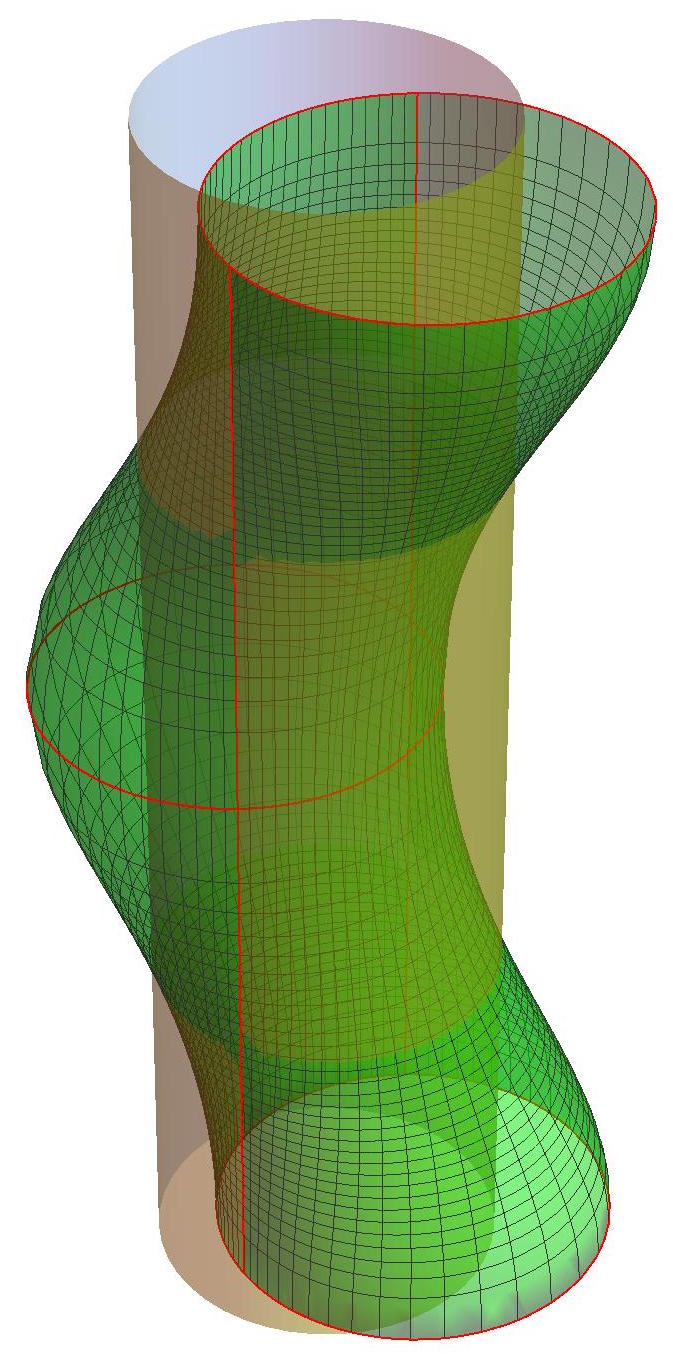}
\caption{Tall rectangles can also be seem as the intersection of a 
cylindrical surface of $\RR^3$ with the cylinder $\{(x,y,t) \in \RR^3 \; : \; x^2+y^2<1\}$.}
\label{graph-1}
\end{center}
\end{figure}

Observe that, for a fixed $d \in (0,1)$, the parametrization $\Upsilon_d$ is defined on $(d_1,1/d_1) \times (-1,1)$
and the image of this extension is an annulus in $\R^3$, see Figure \ref{annulus1}. Applying an $180^o$-degree rotation around the line
(geodesic) passing through the points $(-{\rm i},h_d/2)$ and $({\rm i},h_d/2)$,
then we get another annulus of height $2 h_d$ which is the fundamental
piece of a singly periodic surface in $\R^3$. The  intersection of this periodic
surface with the solid cylinder $\{ (x,y,t) \in \RR^3 \; : \; x^2+y^2 < 1\}$
consists of an infinite number copies of the tall rectangle $\Sigma_d$, see Figure \ref{graph-1}.

\section{Jacobi fields on parabolic catenoids}
In the remaining sections of this paper we initiate the study of the Jacobi operator on a parabolic catenoid.
This is a necessary step before studying the broader space of minimal disks with a similar asymptotic structure.
In this section we consider the Jacobi fields on $\fD$ which are generated by deformations into catenoids or
tall rectangles. The next section is a brief introduction to some more general aspects of the analysis of the Jacobi operator.

\subsection{Deformations to catenoids} We first analyze the deformations of $\fD$ into catenoids $\calC_{q, h}$. A
thorough analysis of the Jacobi operator and local deformation theory for catenoids $\calC_\kappa$ appears in
our earlier paper \cite{fmmr}, but we consider here this limiting case and find an explicit expression for the Jacobi field
on the parabolic catenoid arising from the `regeneration' of this surface to the degenerating sequence of catenoids.  This is
a bit complicated because the rotational symmetry of the catenoids is lost in the limit. For this reason, we
do all computations in the rectangular coordinates $(x,y,t)$ on $\hs \times \R$. 

Recall from Section \ref{cat} the conformal parametrization
$$
X_\ka(\theta,t)=\left( r_\ka(t) {\rm e}^{{\rm i} \; \sqrt{\ka} \;\theta}, t\right), \quad \ka>0,
$$
where %$r_\ka(t)$ is a solution of the ODE
\begin{equation} \label{eq:radius}
\ka= \frac{(1-r_\ka(t)^2)^2}{4 r_\ka(t)^2}+\left( \frac{r_\ka'(t)}{r_\ka(t)}\right)^2.
\end{equation}
We choose the solution for which $r < 1$ when $t \in \left(-\frac{h_\ka}{2},\frac{h_\ka}{2} \right)$. 
We also know that 
$$
\min r_\ka := r_0 =r_\ka(0)=\sqrt{\ka+1}-\sqrt{\ka}\leq r_\ka(t) < 1.
$$

\begin{figure}[htbp]
\begin{center}
\includegraphics[width=8.5cm]{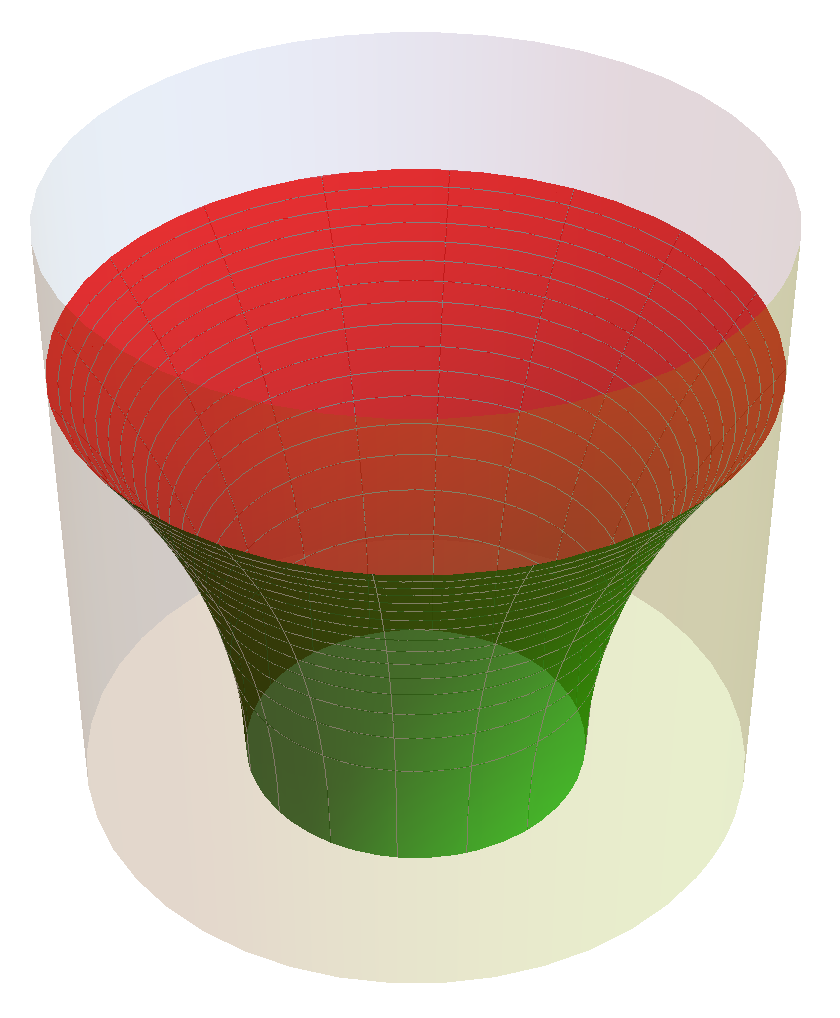}
\caption{Half a catenoid is a graph over the annulus $A(r_0,1)$.}
\label{graph}
\end{center}
\end{figure}

Now write the catenoid as a bigraph over the planar annulus $A(r_0,1)$, see Figure \ref{graph},
with $t$ a function of $r$.  Then
$$
\ka= \frac{(1-r^2)^2}{4 r^2}+\left( \frac{1}{r \; t'(r)}\right)^2,
$$
or equivalently, $ t'(r)=2 (4 \ka r^2-(1-r^2)^2)^{-1/2}$, whence 
$$
t_\ka(r)= \int_{r_0}^r \frac{2 \; du}{\; \sqrt{4 \ka u^2-(1-u^2)^2} \;}.
$$

Before differentiating this with respect to $\ka$, it is convenient to set $u=(r-r_0) s+r_0$, which
yields 
$$
t_\ka(r)=2 (r-r_0) \int_{0}^1 \frac{ds}{\; \sqrt{4 \ka ((r-r_0) s+r_0)^2-(1-((r-r_0) s+r_0)^2)^2} \;}.
$$

We next transform to the upper half-space times $\R$, using the conformal diffeomorphism $g:\D \rightarrow \hs$
from the preliminaries. Write
$$
\mu_0:=-{\rm i} \; g(r_0)=\frac{1+\sqrt{\ka}-\sqrt{1+\ka}}{1-\sqrt{\ka}+\sqrt{1+\ka}} \in (0,1). 
$$
We obtain a representation of $\mathcal{C}_\ka$ as a bigraph over the complement of a disk in the half-plane
$$
\Omega_\ka:=\hs \setminus D({\rm i} \, (\mu_0+\mu_0^{-1})/2,(\mu_0^{-1}-\mu_0)/2) \subset \hs,
$$
where $D(z_0,r):= \{ z \in \C \, : \, |z-z_0|<r\}.$  

Next apply the horizontal dilation $T_{1/\mu_0}$ on $\hs$ which carries $z \mapsto \frac{1}{\mu_0} \,z$,  \label{T}and write
$\whC_\ka:= \left(T_{1/\mu_0} \times {\rm Id}_\R\right)(C_\ka).$  As described earlier, the parabolic  catenoid
is the limit of the family $\whC_\ka$, as $\ka \to 0$.  Note that
$$
\lim_{\ka \to 0} T_{1/\mu_0} (\Omega_\ka) = \widetilde{M} =\{ z = x+{\rm i} y \in \hs \; : \; 0<y \leq 1\}.
$$
Given a point $z \in  \widetilde{M}$, set $r= \left\| g^{-1}(T_{\mu_0}(z)) \right\|$; 
the corresponding point on the dilated catenoid is $(x,y,t_\ka(r))$. % where
%$$r=\frac{\sqrt{2 \left(\left(1-\sqrt{\ka+1}\right)
%   \left(x^2+y^2\right)^2+\sqrt{\ka+1}+1\right)+\ka \left(2 \left(x^2-1\right)
%   y^2+\left(x^2+1\right)^2+y^4\right)}}{\left(\sqrt{\ka+1}-1\right)
%   x^2+\left(\sqrt{\ka+1}-1\right) y^2+2 \sqrt{\ka} y+\sqrt{\ka+1}+1}$$
This gives the family of minimal immersions 
$$
\Phi_\ka(x,y)=(x,y,t_\ka(x,y)).
$$
To compute $\left. d \Phi_\ka/d \ka \right|_{\ka=0}$, it is first necessary to compute  $\left. d t_\ka/d \ka \right|_{\ka=0}$.
Write the integral formula for $t_\ka$ above as $t_\ka(x,y)=\int_0^1 G(x,y,s,\ka) ds$ and expand $G(x,y,s,\ka)$ in
powers of $\sqrt{\ka}$: 
\begin{multline*}
G(x,y,s,\ka)=a_0(x,y,s)+ a_1(x,y,s) \sqrt{\ka}+ a_2(x,y,s) \ka +o(\ka^{3/2}). 
\end{multline*}
These first few coefficients are given by
$$
a_0(x,y,s)= \frac{1-y}{\sqrt{s(1-y)(2-s+s y)}},
$$
\begin{multline*}
a_1(x,y,s) \\ =\frac{-s^2 y^3+3 s^2 y^2-3 s^2 y+s^2-3 s y^2+6 s y-3 s+y^2-2 y+1}{2 (s y-s+2)
  \sqrt{s   (1-y) (s y-s+2)}}
\end{multline*}
and
\begin{multline*}
  a_2(x,y,s) \\ =\frac{-2 s^4 (y-1)^5+2 s^3 (y-5) (y-1)^4+s^2 (10 y-13) (y-1)^3}{8 (s
   y-s+2)^2 \sqrt{s (1-y) (s y-s+2)}} \\
  +\frac{ 2 s (y-1) \left(y \left(x^2+8
   y-6\right)-2\right)+y \left(4 x^2+(5-3 y) y+7\right)-9}{8 (s
   y-s+2)^2 \sqrt{s (1-y) (s y-s+2)}}
   \end{multline*}
A straightforward but increasingly tedious computation gives 
\begin{eqnarray}
\int_0^1 a_0(x,y,s) ds &= & \arccos y, \\
\int_0^1 a_1(x,y,s) ds &=& 0, \\
\int_0^1a_2(x,y,s) ds &= &\frac 14 \left( \frac{x^2 y}{\sqrt{1-y^2}}-\arccos y\right).
\end{eqnarray}
These yield
$$
\lim_{\ka \to 0} \Phi_\ka(x,y)=\Phi_0(x,y)=(x,y,\arccos y),
$$
and
$$
\left.\frac{d \Phi_\ka}{d \ka} \right|_{\ka=0}(x,y)=\left(0,0,\frac 14 \left( \frac{x^2 y}{\sqrt{1-y^2}}-\arccos y\right) \right).
$$
   
   Finally, if we parametrize half of the parabolic  catenoid as $\Phi_0(x,y)=(x,y,\arccos y)$, then the Gauss map 
   is given by $\nu(x,y)=(0,-y^2,-\sqrt{1-y^2}).$ Hence, the Jacobi field that we are looking for is
   $$w(x,y)= \nu \cdot \left( \left.\frac{d \Phi_\ka}{d \ka} \right|_{\ka=0}\right)=\frac 1 4 (\sqrt{1-y^2} \arccos y-x^2 y).$$
   Using the parametrization $\widehat F: \R \times \left(-\frac \pi2, \frac \pi2 \right) \rightarrow \h^2 \times \R$, 
   $$
   \widehat F(x,t)=(x,\cos t,t),
   $$
then this Jacobi field equals $w(x,t)=\frac 14 (t \sin t-x^2 \cos t).$
\begin{remark}
The parabolic  catenoid in \eqref{daniel} and the one above differ by the vertical translation
$t \mapsto t-\pi/2$.  The Jacobi field on \eqref{daniel}  is
$$
w(x,t)=\frac{1}{8} \left((\pi -2 t) \cos (t)-2 x^2 \sin (t)\right).
$$
\end{remark}

\subsection{Deformations to tall rectangles}
We next compute the Jacobi field on a parabolic catenoid associated to the variation of this surface
into the family of tall rectangles. More specifically, we compute the variation associated to the family $Y_d$ 
defined at the end of section 5 which converge to a given parabolic  catenoid $\fD$. 

% The following are the possibilities for the geometric limits of this family of disks: 
 
%  \begin{enumerate}[a.]
% \item Fixing the height $h$, then:
%  \begin{enumerate}[{a}.1]
%  \item $\Sigma_{\sigma h} $ converges to the vertical segment 
%  $\{q_1\}  \times \left[-\frac h2, \frac h2  \right]$ as the arc $\sigma$
%  converges to the point $\{ q_1\}.$
%  \item $\Sigma_{\sigma h} $ converges to the horizontal slices
%  $\h^2 \times \{-h/2, h/2\} $ joint by the vertical segment 
%  $\{q_1\} \times \left[-\frac h2, \frac h2  \right],$ as $\sigma$ 
%  converges to the whole $\partial_\infty \h^2.$
%  \end{enumerate}
%  \item Fixing the arc $\sigma$, then:
%   \begin{enumerate}[{b}.1]
%    \item $\Sigma_{\sigma h} $ converges to the totally geodesic
%    plane $[q_1,q_2] \times \R$, as $h \to \infty.$
%    \item $\Sigma_{\sigma h} $ diverges to the region in 
%    $\partial_\infty(\h^2 \times \R)$ bounded by 
%    $$ \left( \sigma \times \left\{ -\frac \pi2, \frac \pi2 \right\}\right) 
%    \cup  \left( \left\{  q_1,q_2 \right\}\times \left[ -\frac \pi2, \frac \pi2 \right] \right), $$
%    as $h \to \pi.$
%    \end{enumerate}
%  \end{enumerate}

Recall the parametrization \eqref{Upsilon} for $\Sigma_d$. 
% The parametrization of ll cathe minimal disk $\Sigma^+_d$ is:
% $$\Upsilon_d^+ : (d_1,1) \times (-1,1) \rightarrow \h^2 \times \R$$
% %%
% $$ \Upsilon_d^+(x, y)=
% \left(\frac{x-x y^2}{x^2 y^2+1},\frac{\left(x^2+1\right) y}{x^2 y^2+1}, \lambda_d(x)\right)$$
The top and bottom halves of $\Sigma_d$ are graphs over the lunette $L_d \subset \D$
between the circular arc $\gamma_d$ passing through $\pm i$ and $d_1$, and the boundary arc
$\gamma_1$ passing through $\pm i$, $1$. Set $\widehat L_d=g(L_d)$; this is the region in the half-plane $\hs$
between the segment $[-1,1]$ and the circular arc passing through $\pm 1$ and $\mu_1$, where 
$$
\mu_1= {\rm i} \; \frac{\sqrt{1+d}-\sqrt{1-d}}{\sqrt{1+d}+\sqrt{1-d}}.
$$
Defining $T_s$ as before, then the limit of the domains $T_{1/\mu_1}(\widehat L_d)$ as $d \to 0$ is the strip $M$
where $0 < y < 1$. 

Given $z = (x,y) \in M$, $z \in T_{1/\mu_1}(\widehat L_d)$ for all $d<d'$ if $d'$ is small enough.
% $T_{\mu_1}(x,y)$.
The transformation of $\Sigma_d$ is parametrized by the function
$$
G_d(x,y)=(x,y,\lambda_d(X)),
$$
where $X=X(x,y,d)$ is the positive solution to the quadratic system:
$$
\begin{cases} \frac{X(1-Y^2)}{1+X^2 Y^2}= \mu_1 x \\
\mbox{} \\ 
\frac{Y(1+X^2)}{1+X^2 Y^2}= \mu_1 y \end{cases}
$$
Now recall the formula \eqref{lambdad} for $\lambda_d(X)$, and change variables in it by setting $r=(X-d_1) s+d_1$. This gives
\[
  \lambda_d(X)  =  \int_0^1 H(x,y,s,d) ds.
\]
where
\begin{multline*}
H(x,y,s,d) = \frac{2 \; (X-d_1)ds}{\sqrt{d^2(((X-d_1) s+d_1)^2+1)^2-(((X-d_1) s+d_1)^2-1)^2}} 
\end{multline*}
Now proceed as in the previous subsection by expanding %$H(x,y,s,\sqrt{d})$ at $d=0$:
$$
H(x,y,s,\sqrt{d})= h_0(x,y,s)+h_1(x,y,s) \sqrt{d}+h_2(x,y,s) d+ o(d^{3/2}),
$$
where
$$
h_0(x,y,s)= \frac{1-y}{\sqrt{s(1-y)(2-s+s y)}},
$$
$$
h_1(x,y,s)=\frac{(1-y)^{3/2} \left(s^2 y-s^2+3 s-1\right)}{2 \sqrt{s} (s y-s+2)^{3/2}},
$$
and
\begin{multline*}
  h_2(x,y,s)= \\
  \frac{-2 s^4 (y-1)^5+2 s^3 (y-5) (y-1)^4+s^2 (10 y-17) (y-1)^3}{8 (s    (y-1)+2)^2 \sqrt{s (1-y) (s (y-1)+2)}}+ \\
\frac{2 s (y-1) \left(-\left(x^2+14\right) y+8 y^2+6\right)-y \left(4 x^2+y (3 y-5)+9\right)+7}{8 (s
  (y-1)+2)^2 \sqrt{s (1-y) (s (y-1)+2)}}.
\end{multline*}
   
By a straightforward computation,
\begin{multline*}
  \int_0^1 h_1(x,y,s) ds=0, \quad \mbox{and} \\
  \int_0^1 h_2(x,y,s) ds=\frac{1}{4} \left( \arccos(y)-\frac{ x^2
      y}{\sqrt{1-y^2}}\right).
  \end{multline*}
  Therefore, just as at the end of the previous subsection, the corresponding Jacobi field is
  $$
  \widehat w(x,y)=-\frac 1 4 (\sqrt{1-y^2} \arccos y-x^2 y).
  $$
  Note the unexpected fact that this is equal to $-w(x,y)$, where $w$ is the Jacobi field associated to
  the deformation of $\fD$ into catenoids.

\section{The Jacobi operator on parabolic catenoids}
We finally turn to the task of describing the beginnings of the analytic theory of the Jacobi operator on the
parabolic catenoid $\fD _{\infty,1}$ given by parametrization \eqref{daniel} for $\lambda=1$. 

\noindent{\bf Coordinate vector fields and metric:}   Via the parametrization $\Psi$, the coordinates
$(x,t)$ on $\Sigma$ induce the coordinate vector fields
\[
\Psi_*(\del_x) = X_1 = (1, 0, 0), \qquad \Psi_*(\del_t) = X_2 = ( 0, \cos t, 1).
\]
The metric coefficients are
$$
g_{11} = X_1 \cdot X_1 = 1/\sin^2 t, \quad g_{12} = X_1 \cdot X_2 = 0, $$
$$\mbox{and}\quad g_{22} = X_2 \cdot X_2 = (\cos^2 t/\sin^2 t) +1= 1/\sin^2 t,
$$
(thus displaying the conformality of $\Psi$).

\medskip

\noindent{\bf Jacobi operator:} 
The unit normal to $\Sigma$ at $\Psi(x,t)$ is
\[
\nu(x,t) = (0, \sin^2 t, -\cos t),
\]
whence $\mathrm{Ric}_{\h^2 \times \R} (\nu, \nu) = -g_{\h^2}( (0, \sin^2 t), (0, \sin^2 t) ) = - \sin^2 t$.  To compute the
Jacobi operator we can avoid computing $|A|^2$ directly by the following observation.

We first compute the Jacobi field corresponding to the family $\lambda \mapsto \Sigma_{\infty,0,\lambda}$.   Indeed,
\[
\left. \frac{d\,}{d\lambda}\right|_{\lambda = 1} \Psi_\lambda(x,t) = (x, \sin t, 0),
\]
so the normal component of this, which is the Jacobi field we seek, equals
\[
\psi = \nu \cdot (x, \sin t, 0) = \sin t.
\]
This vanishes simply at both $t=0$ and $t=\pi$ and is $L^2$ on any portion $|x| \leq C$ since the measure
equals ${dx dt}/ \sin^2 t$, but is not $L^2$ on $\Sigma$. Now, $\Delta_g = \sin^2 t ( \del_x^2 + \del_t^2)$ and 
the Jacobi operator equals 
\[
L = - \left( \Delta_g + \mathrm{Ric}(\nu, \nu) + |A|^2\right).
\]
Writing out the equality $L \psi = 0$ with $\psi= \sin t$, we obtain
\[
\sin^2 t (- \sin t) + (- \sin^2 t + |A|^2) \sin t = 0 \Rightarrow |A|^2 = 2 \; \sin^2 t
\]
This shows that
\[
L = -\sin^2 t ( \del_x^2 + \del_t^2 + 1).
\]
This is a nonnegative operator: indeed, $L = \sin^2 t L_0$, where $L_0 = -\del_x^2 - \del_t^2 - 1$, and its action on
$L^2( (\sin t)^{-2} dxdt)$ is equivalent to the action of $L_0$ on $L^2( dxdt)$ with Dirichlet boundary conditions at
$t=0, \pi$. This latter operator is nonnegative since $-\del_t^2 - 1$ with Dirichlet conditions is nonnegative on $[0,\pi]$.
The fact that the spectrum of $L$ lies in $\R^+$ also follows from the existence of the nonnegative solution $\psi$
to $L\psi = 0$.

The function $\tilde{u}(x,t) = x \sin t$ is another non-$L^2$ solution to $L\tilde{u} = 0$. It is not hard to show
that $u$ and $\tilde{u}$ span the full space of tempered solutions to the Jacobi equation which vanish at
$t = 0, \pi$. 

\medskip

\noindent{\bf Mapping properties of the Jacobi operator:}  We next describe some aspects of the
mapping properties of $L$ on the infinite strip $S = \RR_x \times [0,\pi]_t$. The remarks here
are meant to be preparatory to a deeper study of the deformation theory of tall rectangles, to which
we shall return in a work in progress.   We consider two main questions:
\begin{itemize}
\item[i)] Find classes of functions $\phi_\pm(x)$ such that the problem
$Lu = 0$ on $S$, $u(x,\pi) = \phi_+(x)$, $u(x,0) = \phi_-(x)$ 
is solvable;
\item[ii)] Find a class of functions $f$ on $S$ for which we can solve $Lu = f$ with $u = 0$ at $t=0, \pi$, and $u \to 0$
  as $x \to \pm \infty$.
  \end{itemize}

We analyze these questions using the Fourier transform in $x$. Writing the dual variable as $\xi$, then question i)
leads to the study of the two families of problems
\begin{equation}
  \begin{split}
    & \hat{L}_\xi \hat{u} (\xi,t) :=  \sin^2 t ( -\del_t^2 + \xi^2 - 1) \hat{u} = 0, \\
    & \hat{u}(\xi, \pi) = \hat{\phi}_+(\xi),\ \hat{u}(\xi,0) = \hat{\phi}_-(\xi)
  \end{split}
  \label{ihdp}
\end{equation}
and
\begin{equation}
  \hat{L}_\xi \hat{u}(\xi, t) = \hat{f}(\xi, t), \qquad \hat{u}(\xi, 0) = \hat{u}(\xi, \pi) = 0.
  \label{ihe}
  \end{equation}
  
Consider \eqref{ihdp} first.  For $\xi \neq 0$, there exist two functions $v_\pm(\xi, t)$ which satisfy $\hat{L}_\xi  v_\pm(\xi, t) = 0$ and
\[
v_+(\xi,0) = 0,\ v_+(\xi,\pi) = 1, \quad v_-(\xi,0) = 1,\ v_-(\xi,\pi) = 0,
\]
namely
\begin{equation*}
\begin{cases}
  & v_+(\xi,t) = \sin ( (1-\xi^2)^{1/2} t )/ \sin( (1-\xi^2)^{1/2} \pi), \\
  & v_-(\xi,t) = \sin ( (1-\xi^2)^{1/2} (\pi-t) )/ \sin( (1-\xi^2)^{1/2} \pi)
\end{cases}
\end{equation*}
when $|\xi|< 1$,
\[
  v_+(\pm 1,t) =  t/\pi,\ v_-(\pm 1, t) = 1 - t/\pi,
\]
and
\begin{equation*}
\begin{cases}
& v_+(\xi,t) = \sinh ( (\xi^2-1)^{1/2} t)/\sinh ( (\xi^2-1)^{1/2} \pi), \\
&  v_-(\xi,t) = \sinh ( (\xi^2-1)^{1/2} (\pi-t))/\sinh ( (\xi^2-1)^{1/2} \pi)
\end{cases}
\end{equation*}
for $|\xi| > 1$.   These functions are clearly holomorphic when $\xi \in \mathbb C \setminus \{0, \pm 1\}$, but
although they appear to be branched at $\xi = \pm 1$, they are single-valued at these points so are holomorphic
away from $\xi = 0$, where they have a double pole. 

Now suppose, for example, that $\hat{\phi}_\pm(\xi)$ are functions in $L^1$ such that $\phi_\pm(\xi) \xi^{-2} \in L^1_{\mathrm{loc}}$. 
Then the solution to problem i) is
\[
u(x,t) = \mathcal F^{-1} \left( \hat{\phi}_+(\xi) v_+(\xi,t) + \hat{\phi}_-(\xi) v_-(\xi,t) \right),
\]
where $\mathcal{F}$ is the Fourier transform.   It is not hard to check that under these hypotheses, $u(x, t)$ is continuous
on $\mathbb R \times [0,\pi]$,  $u(x,t) \to 0$ uniformly in $t \in [0,\pi]$ as $x \to \pm \infty$, and furthermore, that
\[
\int_{\RR}  u(x,\pi)\, dx = \int_{\RR} u(x, 0)\, dx = 0.
\]
It is clearly possible to choose functions $\phi_\pm$ satisfying these constraints but so that $u(x,\pi) \neq u(x,0)$ for
every $x \in \RR$.  This implies that there are infinitesimal deformations where the difference of the heights of the
two boundary curves may vary, though the average of the difference of their heights equals $\pi$.

There are also some interesting constraints on Jacobi fields. Indeed, let
$$
S_r=\{ \Psi(x,t) \; : \; (x,t) \in [-r,r] \times [0,\pi] \}
$$
denote the truncated surface.  Consider the basic Jacobi field $u(x,t) = \sin t$, and suppose that $w$ is any other Jacobi field, i.e.,
$Lw = 0$, which has sufficient decay as $|x| \to \infty$ for the following computations to make sense (such Jacobi fields certainly
exist by virtue of the preceding calculations.) We then compute that  
\begin{multline*}
0 = \int_{S_r}  \left( (Lu) w - u (Lw)\right) \, \frac{1}{\sin^2 t} dxdt = \\  \qquad \int_{-r}^r \left( w(x,\pi) + w(x,0) \right) \, dx +
\int_0^\pi \sin t \,\left(w_x(r,t)-w_x(-r,t)\right) dt. \label{eq:A}
\end{multline*}
We conclude that Jacobi fields (at least the well-behaved ones) must satisfy the `moment condition'
$$
\int_{-\infty}^\infty (w(x,\pi)+w(x,0)) dx+ \int_0^\pi \sin t \,(w^+_x(t)-w^-_x(t)) dt=0,
$$
where $\displaystyle w^\pm_x(t):= \lim_{x \to \pm \infty} w_x(x,t).$    The precise geometric meaning of this
is not evident. 

To convert these infinitesimal statements into statements about minimal surfaces near to $\Sigma_{q, \tau, \lambda}$,
it is necessary to solve the inhomogeneous problem ii).   The details of this proceed in an unsurprising fashion: passing to
the Fourier transform again, there is a Green function $\hat{G}(\xi, t, t')$ for $\hat{L}_\xi$ for $\xi \neq 0$,
and this can be used to solve $Lu = f$ for a broad collection of functions $f$.   By this linear theory and standard
use of the implicit function theorem, the Jacobi fields discussed earlier can be integrated to nearby minimal surfaces.
 
The point of all of this is the following.  There exist deformations of $\Sigma_{q, \tau, \lambda}$ which deform the top
and bottom boundary curves to $\Gamma_\pm$, but which fix the vertical line connecting them. Although
it might be natural to conjecture that any such deformation has top and bottom boundary separated at exactly
distance $\pi$, i.e., $\Gamma_+(x) - \Gamma_-(x) = \pi$ for all $x$, we have shown 
that this is not true.

\end{document}